\documentstyle[12pt,twoside]{amsart}

\title
{Smoothness and jet schemes}
\author{Shihoko Ishii} 

\address{Department of Mathematics, Tokyo Institute of
Technology, Oh-Okayama, Meguro, 152-8551 Tokyo, Japan
\newline
e-mail : ishii.s.ac@@m.titech.ac.jp}

\newcommand{\bZ}{{\Bbb Z}}

\newcommand{\bN}{{\Bbb N}}
\newcommand{\AAA}{{\widehat{{\Bbb A}_{k}^N}}}
\newcommand{\bA}{{\Bbb A}}

\newcommand{\bx}{{\bf x}}

\newcommand{\Spec}{\operatorname{Spec}}

\newcommand{\Hom}{\operatorname{Hom}}
\newcommand{\ord}{\operatorname{ord}}

\newcommand{\stm}{{\Spec k[t]/(t^{m+1})}}

\let \cedilla =\c
\renewcommand{\c}[0]{{\mathbb C}}  

\renewcommand{\o}[0]{{\mathcal O}} 

\newcommand{\spec}[0]{\operatorname{Spec}}

\newcommand{\ed}{{\operatorname{embdim}}}

\def\to {\longrightarrow}

\newtheorem{thm}{Theorem}[section]

\newtheorem{lem}[thm]{Lemma}
\newtheorem{cor}[thm]{Corollary}
\newtheorem{prop}[thm]{Proposition}

\theoremstyle{definition}
\newtheorem{defn}[thm]{Definition}

\newtheorem{say}[thm]{}
\newtheorem{exmp}[thm]{Example}

\newtheorem{rem}[thm]{Remark}

\theoremstyle{remark}

\begin{document}
\maketitle
\footnote{partially supported by Grant of JSPS}

\begin{abstract}
This paper shows some criteria for a scheme of finite type over an 
algebraically closed field to be non-singular in terms of jet 
schemes. 
 For the base field of  characteristic zero,  
 the scheme is non-singular if and only if one of the  truncation morphisms 
 of its jet schemes 
 is flat.
  For the positive characteristic case, 
  we obtain a similar characterization under the reducedness condition on the 
  scheme.
  We also obtain by a simple discussion that  the scheme is non-singular if and only if one of its jet schemes is 
non-singular. 
\end{abstract}

\section{Introduction}
\noindent
In 1968 John F. Nash introduced the  jet schemes and the arc space of 
an algebraic and an analytic variety  and  
 posed the Nash problem (\cite{nash}).

  The jet schemes and the arc space  are  considered to be something to reflect 
  the nature of the
  singularities of a variety. 
  (The Nash problem itself concerns a connection  between the  arc 
  space and the singularities.)
  By looking at the jet schemes over a variety, 
  we can see  some properties of the 
  singularities of the variety (see \cite{ein}, \cite{e-Mus}, \cite{must01}, 
  \cite{must02}) : for example, if \( X \) is  
  locally a complete intersection variety, the singularities of \( X 
  \) are canonical (resp. 
  terminal) if and only 
  if the jet scheme \( X_{m} \) is irreducible (resp. normal) for every 
  \( m\in \bN \).
  
  For  a non-singular variety \( X \), the jet 
  schemes are distinguished: the \( m \)-jet scheme \( X_{m} \) is 
  non-singular for every \( m\in \bN \) and every truncation morphism \( 
  \psi_{m',m}:
  X_{m'}\to X_{m}\) is smooth with the fiber \( \bA_{k}^{(m'-m)\dim X} \)
   for \( m'>m\geq 0\). 
  Then, it is natural to ask whether  these properties 
  characterize the smoothness of the variety \( X \).
  
  Our results are rather stronger, i.e., only one jet scheme or  one 
  truncation morphism is sufficient to characterize the smoothness of the variety \( X \).
  In this paper we prove the following:
  
\begin{prop}
\label{sm}
Let \( k \) be a field of arbitrary 
characteristic and \( f:X\to Y \) a morphism of \( k \)-schemes.
Then the following are equivalent:
\begin{enumerate}
   \item[(i)] \( f \) is smooth (resp. unramified, \'etale);
   \item[(ii)] For every \( m\in \bN \), the morphism \( f_{m}:X_{m}\to 
    Y_{m} \) induced from $f$ is smooth (resp. unramified, \'etale);
   \item[(iii)] There is an integer \( m\in \bN \) such that 
   the morphism \( f_{m}:X_{m}\to 
    Y_{m} \) is smooth (resp. unramified, \'etale).
\end{enumerate}
\end{prop}

As a corollary of this proposition, we obtain the following:
  
\begin{cor}
\label{smooth}
Let \( k \) be a field of arbitrary 
characteristic.
A scheme \( X  \) of finite type over  \( k \)  is 
smooth if and only if there is \( m\in \bZ_{\geq 0} \) such that \( X_{m} 
\) is smooth.
\end{cor}

\begin{thm}
\label{flat}
Let \( k \) be an algebraically closed field of characteristic zero. 
A scheme \( X \)  of finite type over \( k \) 
  is non-singular if and only if there is  a pair  of 
integers \( 0\leq m<m' \) such that the truncation morphism \( \psi_{m', m}: X_{m'}\to 
X_{m} \) is a flat morphism.
\end{thm}

Here, we note that the assumption of the characteristic of the base 
field in Theorem \ref{flat} is necessary.
We will see a counter example of this statement in positive 
characteristic(Example \ref{ex}).

If we assume that the scheme $X$ is reduced, 
then we have a similar criterion as Theorem \ref{flat}  also for the positive characteristic case.

\begin{thm}
\label{positive}
Let \( k \) be an algebraically closed field of arbitrary characteristic.
Assume the scheme \( X  \)  of finite type over \( k \) is reduced. 
Then $X$ is non-singular   if and only if there is  a pair  of 
integers \( 0< m<m' \) such that 
the 
truncation morphism \( 
\psi_{m',m}:X_{m'}\to X_{m} \) is flat .
\end{thm}

%

This paper is  motivated by  Kei-ichi Watanabe's question. 
The author expresses her hearty thanks to him. 
The author is also grateful to Mircea Musta\cedilla{t}\v{a} for his 
helpful comments and stimulating discussions.

\section{Preliminaries on jet schemes}

In this paper, a $k$-scheme is always a separated scheme  over a field $k$.

\begin{defn}
  Let \( X \) be a scheme of finite type over \( k \)
and $K\supset k$ a field extension.
  A  morphism \( \Spec K[t]/(t^{m+1})\to X \) is called an \( m \)-jet
  of \( X \).
\end{defn}

\begin{say}
\label{field}
  Let \( X \) be a scheme of finite type over \( k \).
  Let \( {\cal S}ch/k \) be the category of \( k \)-schemes  
   and \( {\cal S}et \) the category of sets.
  Define a contravariant functor  \( {\cal F}_{m}^X: {\cal S}ch/k \to {\cal S}et \)
  by 
$$
 {\cal F}_{m}^X(Y)=\Hom _{k}(Y\times_{\Spec k}\Spec k[t]/(t^{m+1}), X).
$$
  Then, \( {\cal F}_{m}^X \) is representable by a scheme \( X_{m} \) of finite
  type over \( k \), that is
$$
 \Hom _{k}(Y, X_{m})\simeq\Hom _{k}(Y\times_{\Spec k}
\Spec k[t]/(t^{m+1}), X).
$$ 
   This \( X_{m} \) is called the {\it scheme of \( m \)-jets} of \( X 
   \) or the {\it \( m \)-jet scheme} of \( X \).
  For \( m<m' \) the canonical surjection \( k[t]/(t^{m'+1})\to k[t]/(t^{m+1}) \)
  induces a morphism \( \psi^X_{m',m}:X_{m'}\to X_{m} \), 
  which we call a truncation morphism.
  In particular, for \( m=0 \)  \( \psi^X_{m,0}:
  X_{m}\to X \) is denoted by \( \pi^X_{m} \).
  We denote   \(  \psi^X_{m',m} \) and \( \pi^X_{m} \) by  \(  
  \psi_{m',m} \) and \( \pi_{m} \), respectively, if there is no risk 
  of confusion.
  By \ref{field}, a point \( z \in X_{m} \)  gives an \( m \)-jet \( 
  \alpha_{z}:
  \Spec K[t]/(t^{m+1})
  \to X \) and \( \pi^X_{m}(z)=\alpha_{z}(0) \),
  where \( K \) is the residue field at \( z \) and \( 0 \) is 
  the point of \( \spec K[t]/(t^{m+1})  \). 
  From now on we denote a point \( z \) of \( X_{m} \) and the 
  corresponding \( m \)-jet \( \alpha_{z} \) by the common symbol \( 
  \alpha \). 
\end{say}  

\begin{say}
  The canonical inclusion \(  k\to k[t]/(t^{m+1}) \)  induces a 
  section \( \sigma^X_{m}:X \hookrightarrow X_{m}  \) of \( \pi^X_{m} \).
  The image \( \sigma^X_{m}(x) \) of a point \( x\in X \) is the 
  trivial \( m \)-jet at \( x \) and is denoted by \( x_{m} \).
\end{say}

\begin{say}
 Let \( f:X\to Y \) be a morphism of \( k \)-schemes.
  Then the canonical morphism \( f_{m}:X_{m}\to Y_{m} \) 
  is induced for every \( m\in \bN \) such that the 
  following diagram is commutative:
  \[ \begin{array}{ccc}
      X_{m}& \stackrel{f_{m}}\longrightarrow & Y_{m}\\
      \pi^X_{m} \downarrow\ \ \ \ \ & & \ \ \ \ \downarrow \pi^Y_{m}\\
      X & \stackrel{f}\longrightarrow & Y\\
      \end{array}. \]
      Pointwise, for \( \alpha\in X_{m} \) , \( f_{m}(\alpha)\) is the $m$-jet  \[ f\circ \alpha:
      \spec K[t]/(t^{m+1})\stackrel{\alpha}\longrightarrow X \stackrel{f}\longrightarrow Y. \]
\end{say}

\section{Proof of Proposition \ref{sm}}

\noindent
[{\it Proof of Proposition \ref{sm}}]
  (i)\( \Rightarrow \) (ii):  This implication for smooth and \'etale cases is already mentioned in \cite{BL}  and  \cite{EM}.
  For the reader's convenience, the proof is included here.
    Assume for an integer \( m\geq 0 \), a commutative diagram of \( k \)-schemes:
\[ \begin{array}{ccc}
X_{m}& \stackrel{f_{m}}\longrightarrow & Y_{m}\\
\uparrow& & \uparrow\\
Z' &\hookrightarrow& Z\\
\end{array}  
  \] 
  is given, where \( Z'\hookrightarrow Z \) is a closed immersion of affine schemes 
  whose defining ideal is nilpotent.
  This diagram is equivalent to the following commutative diagram:
  \[ \begin{array}{ccc}
X& \stackrel{f}\longrightarrow &Y \\
\uparrow& & \uparrow\\
Z'\times \stm &\hookrightarrow& Z\times \stm\\
\end{array}. \]
Here, we note that \( Z'\times \stm \hookrightarrow Z\times \stm \) 
is a closed subscheme with the nilpotent defining ideal.
If \( f \) is smooth (resp. unramified, \'etale), 
there exists a (resp. there exists at most one, there exists a unique) 
morphism \( Z\times \stm \to X \) which makes the two triangles 
commutative.
This is equivalent to the fact that there exists a (resp. there exists at most one, 
there exists a unique)
morphism \( Z \to X_{m} \) which makes the two triangles in the first 
diagram 
commutative.

\noindent
(ii)\( \Rightarrow \) (iii): trivial.

\noindent
(iii)\( \Rightarrow \) (i):
Assume a commutative diagram,
  \begin{equation}\label{d1}
   \begin{array}{ccc}
X& \stackrel{f}\longrightarrow &Y \\
\varphi\uparrow\ \ \ & & \ \ \ \uparrow\psi\\
Z' &\hookrightarrow& Z \\
\end{array} 
\end{equation}
is given, where \( Z'\hookrightarrow Z \) is a closed immersion of affine schemes  whose defining ideal is nilpotent.
  For an integer \( m\geq 0 \), by composing with the sections 
  \(\sigma_{m}^X: X\hookrightarrow X_{m} \), \( \sigma_{m}^Y: 
  Y\hookrightarrow Y_{m} \), we obtain the commutative diagram: 
 \begin{equation}\label{d2} 
 \begin{array}{ccc}
  X_{m}&\stackrel{f_{m}}\longrightarrow& Y_{m}\\
  \cup& & \cup\\
X& \stackrel{f}\longrightarrow &Y \\
\varphi\uparrow\ \ \ & & \ \ \ \uparrow\psi\\
Z' &\hookrightarrow& Z \\
\end{array}. \end{equation}

Now, if \( f_{m} \) is smooth (resp. unramified, \'etale), 
there exists a (resp. exists at most one, exists a unique )
morphism \( Z\to X_{m} \) such that the two triangles are commutative 
in the diagram (\ref{d2}).
By composing this morphism $Z\to X_m$  with $\pi_m^X:X_m\to X$, 
we obtain that there exists a (resp. exists at most one, exists a unique )
morphism \( Z\to X \) such that the two triangles in the lower rectangle are commutative.
\( \Box \)

\vskip.5truecm
\noindent
[{\it Proof of Corollary \ref{smooth}}]
In Proposition \ref{sm}, let \( Y=\spec k \). \( \Box \)

\section{jet schemes of a local  analytic scheme}

For the proofs of the theorems, here we set up the jet schemes for local  analytic schemes.
Let $k$ be an algebraically closed field of arbitrary characteristic. 
The representability of the following functor follows from  \cite{voj}. 
Here, we show the concrete form of the  scheme representing the functor.

\begin{prop}
  Let  $\AAA$ be the affine scheme $\Spec \widehat{\o_{\bA^N,0}}$, where $\o_{\bA^N,0}$ is the local ring of the origin $0\in \bA_k^N$ and  $\widehat{\o_{\bA^N,0}}$ is the completion of $\o_{\bA^N,0}$ at the maximal ideal.
  Let ${\cal F}_m^{\AAA}: Sch/k \to Set$ be the functor from the category of $k$-schemes to the category of sets defined as follows:
  
$${\cal F}_m^{\AAA}(Y):=\Hom_k(Y\times_{\spec k} \stm, \AAA).$$
For a morphism $u:Y\to Z$ in $Sch/k$, 
$${\cal F}_m^{\AAA}(u):\Hom_k(Z\times\stm,\AAA)\to \Hom_k(Y\times\stm,\AAA) $$
is defined by $f\mapsto f\circ (u\times id)$.

Then, ${\cal F}_m^{\AAA}$ is representable by the scheme 
$$(\AAA)_m:=\Spec k[[x_{0,1}, x_{0,2},\ldots,x_{0,N}]][x_{1,1},\ldots,x_{1,N},\ldots,x_{m, 
  1},\ldots,x_{m,N}]$$
   $$  =\spec k[[\bx_{0}]][\bx_{1},\ldots,\bx_{m}], $$
 
  where we denote the multivariables \( ( 
  x_{i,1},x_{i,2},\ldots,x_{i,N})\) by \( \bx_{i} \) for the 
  simplicity of  notation.
\end{prop}

\begin{pf} We may assume that $Y$ is an affine scheme $\spec R$ over $k$.
Then, 
$$\Hom_k(Y\times\stm,\AAA) \simeq \Hom_k(k[[\bx_0]], R[t]/(t^{m+1}))$$
Here we have a bijection:
$$\Hom_k(k[[\bx_0]], R[t]/(t^{m+1}))\simeq \Hom_k(k[[\bx_0]],R)\times R^{mN}$$
by $\varphi\mapsto (\pi_0\circ \varphi, \pi_1\varphi(x_{01}),...,\pi_1\varphi(x_{0,N}),..., \pi_m\varphi(x_{01}),...,\pi_m\varphi(x_{0,N}))$, 
where, $$\pi_i: R[t]/(t^{m+1})\to R\ \ \ \ (i=0,1,...,m) $$  is the projection of $R[t]/(t^{m+1})
=R\oplus Rt \oplus\cdots \oplus Rt^m \simeq R^{m+1}$ 
to the $i$-th factor. 
Indeed it gives a bijection, since we have the inverse map
$$\Hom_k(k[[\bx_0]],R)\times R^{mN}\to \Hom_k(k[[\bx_0]], R[t]/(t^{m+1}))$$
 by 
$$(\varphi_0, a_{1,1},..,a_{1,N},...,a_{m,1},..., a_{m,N})\mapsto \varphi$$
where $\varphi\in \Hom_k(k[[\bx_0]], R[t]/(t^{m+1}))$ is defined as follows:

For $\gamma(x_{0,1}, x_{0,2},\ldots,x_{0,N})\in k[[\bx_0]]$, 
substituting  $\sum_{i=0}^m x_{i,j}t^i$ into $x_{0,j}$ $(j=1,...,N)$ in $\gamma$,
we obtain 
$$\gamma(\sum \bx_i t^i)=\sum_{i=0}^{\infty}\left(\sum_{\sum_{\ell} i_{\ell}=i, 1\leq j_{\ell}\leq N}
\gamma_{i_1,j_1,...,i_s, j_s}x_{i_1,j_1}\cdots x_{i_s,j_s}\right)t^i$$
in $k[[\bx_{0}, \bx_{1},\ldots,\bx_{m}, t]]$, where $\gamma_{i_1,j_1,...,i_s, j_s}\in k[[\bx_0]]$.
Define $\varphi(\gamma)\in R[t]/(t^{m+1})$ by 
$$\varphi(\gamma)=\sum_{i=0}^m\left(\sum_{\sum_{\ell} i_{\ell}=i, 1\leq j_{\ell}\leq N}
\varphi_0(\gamma_{i_1,j_1,...,i_s, j_s})a_{i_1,j_1}\cdots a_{i_s,j_s}\right)t^i.$$

On the other hand,
It is clear that there is a bijection 
$$\Hom_k(k[[\bx_0]][\bx_1,...,\bx_m], R)\simeq \Hom_k(k[[\bx_0]],R)\times R^{mN}$$
by $\varphi\mapsto (\varphi|_{k[[\bx_0]]}, \varphi(x_{1,1}),...,\varphi(x_{1,N}),...,
\varphi(x_{m, 1}),...,\varphi(x_{m,N}))$.
By this, we have
$$\Hom_k(k[[\bx_0]], R[t]/(t^{m+1}))\simeq\Hom_k(k[[\bx_0]][\bx_1,...,\bx_m], R),$$
which implies 
$$\Hom_k(Y\times \spec k[t]/(t^{m+1}), \AAA) \simeq \Hom_k(Y, \spec k[[\bx_0]][\bx_1,...,\bx_m])$$
This completes the proof.
\end{pf}

By this proposition, we have the following:

\begin{cor}
Let $X\subset \AAA$ be a closed subscheme.
Let $I$ be the defining ideal of $X$ in $\AAA$.
Define a functor ${\cal F}_m^X: Sch/k\to Set$ for this $ X$ in the same way  
as in the previous proposition.

For a power series \( f\in k[[\bx_{0}]]\) we define an element 
  \( F_{m}\in k[[\bx_{0}]][ \bx_{1}\ldots,\bx_{m}] \) as follows:
  \[ f(\sum_{i= 0}^m \bx_{i}t^i)=F_{0}+F_{1}t+F_{2}t^2+\cdots+
  F_{m}t^m+\cdots .  \]
  Then, the functor ${\cal F}_m^X$ is represented by a scheme \(  X_{m} \)  defined in 
  \(( \AAA)_{m}= \spec k[[\bx_{0}]][\bx_{1},\ldots,\bx_{m}]\) by 
  the ideal generated by \( F_{i} \)'s \( (i\leq m) \) for all  \( f\in I \).
  (It is sufficient to take \( F_{i} \)'s \( (i\leq m) \) for all generators $f\in I$.)
\end{cor}

\begin{pf}
We use the notation in the proof of the previous proposition.
There, we obtained  bijections :
$$\Hom_k(k[[\bx_0]], R[t]/(t^{m+1}))\stackrel{\Phi}\simeq
 \Hom_k(k[[\bx_0]],R)\times R^{mN}$$
 $$\stackrel{\Psi}\simeq
\Hom_k(k[[\bx_0]][\bx_1,...,\bx_m], R).$$
Here, for $Y=\spec R$, we have the fact that
$${\cal F}_m^X(Y)=\Hom_k(k[[\bx_0]]/I, R[t]/(t^{m+1}))$$
is the subset 
$$\{\varphi:k[[\bx_0]]\to R[t]/(t^{m+1})\mid \varphi(\gamma)=0 \ \ \mbox{for\ generators}\ \gamma\in I\}$$
of $\Hom_k(k[[\bx_0]], R[t]/(t^{m+1}))$.
The condition $\varphi(\gamma)=0$ is equivalent to the conditions $\pi_i\circ\varphi(\gamma)=0$ 
$(i=0,1,...,m)$.
Therefore, this subset is mapped by $\Psi\circ \Phi$ to the subset 
$$\big\{\varphi:k[[\bx_0]][\bx_1,...,\bx_m]\to R\mid \ \varphi(x_{i,j})=a_{i,j}, \ \mbox{for\ generators}\ \gamma\in I, $$
$$\sum_{\sum_{\ell} i_{\ell}=i, 1\leq j_{\ell}\leq N}
\varphi_0(\gamma_{i_1,j_1,...,i_s, j_s})a_{i_1,j_1}\cdots a_{i_s,j_s}=0 \ (i=0,1,...,m)\big\}.$$
Let the ideal $J\subset k[[\bx_0]][\bx_1,...,\bx_m]$ be generated by
$$\sum_{\sum_{\ell} i_{\ell}=i, 1\leq j_{\ell}\leq N}
\gamma_{i_1,j_1,...,i_s, j_s}x_{i_1,j_1}\cdots x_{i_s,j_s}$$ for generators $\gamma\in I,$
then it follows that our subset is equal to 
$$\Hom_k(k[[\bx_0]][\bx_1,...,\bx_m]/J,R).$$

\end{pf}

\begin{rem}
\label{algan}
  Let $X\subset \bA_k^N$ be a closed subscheme  containing the origin $0$, $I_X$  the defining ideal
 and $\widehat X$ the affine scheme $\spec \widehat {\o_{X,0}}$. 
Note that  the defining ideal $I$ of $\widehat{X}$ in ${\AAA}$
 is  generated by $I_X$.  For a polynomial \( f\in k[\bx_{0}] \) we define an element 
  \( F_{m}\in k[\bx_{0}, \bx_{1}\ldots,\bx_{m}] \) in the same way as in the previous corollary.  
 Then \( \widehat X_{m} \) is defined in 
  \(( \AAA)_{m}= \spec k[[\bx_{0}]][\bx_{1},\ldots,\bx_{m}]\) by 
  the ideal generated by \( F_{i} \)'s \( (i\leq m) \) for generators \( f\in I_{X} \).
\end{rem}

\begin{cor}
Under the notation of  Remark\ref{algan}, it follows that 
$$\widehat X_m=\widehat X\times_X X_m.$$
\end{cor}

\begin{pf}
  Note that $F_i\in k[\bx_0, \bx_1,..., \bx_m]$ for a generator $f$ of $I_X$ and $I$ is generated by $I_X$.
 Now the expressions  
  $$X_{m}= \spec k[\bx_{0},\bx_{1},\ldots,\bx_{m}]/ (F_i)_{f\in I_X}$$
  $$\widehat X_{m}= \spec k[[\bx_{0}]][\bx_{1},\ldots,\bx_{m}]/ (F_i)_{f\in I_X}$$
give the required equality.  
\end{pf}

\begin{cor}
\label{equal}
Under the notation of  Remark\ref{algan},
 let $\pi_m^X$ and $\pi_m^{\widehat X}$ be the canonical projections 
 $X_m\to X$ and $\widehat X_m\to \widehat X$, respectively.
 Then, we obtain the  isomorphism of schemes:
 $$(\pi_m^X)^{-1}(0)\simeq (\pi_m^{\widehat X})^{-1}(0).$$
\end{cor}

\begin{cor}
\label{reduction}
  Under the notation of  Remark\ref{algan},
  replacing $X$ by a sufficiently small neighborhood of $0$, we obtain the equivalence that 
  the truncation morphism  $X_{m'}\to X_m$ is flat  
  if and only if  
  the truncation morphism $\widehat X_{m'}\to \widehat X_m$ is  flat.
\end{cor}

\begin{pf}
``Only if" part follows from the base change property for flatness. ``If" part follows from the fact that the homomorphism  
$\o_{X,0}\to \widehat{\o_{X,0}}$ is faithfully flat.
\end{pf}

 \begin{defn}
  \label{weight}
    A monomial \( \bx=\prod_{\ell=1}^d x_{i_{\ell},j_{\ell}}\in k[[\bx_{0}]][\bx_{1},\ldots,\bx_{m}] \)
    is called a monomial of {\it weight} \( w \) if \( w=\sum_{\ell=1}^d 
    i_{\ell}\).
     For an element \( F\in k[[\bx_{0}]][\bx_{1},\ldots,\bx_{m}] \) the 
  order \( \ord F \) is defined as the lowest degree of the monomials 
  in \( \bx_{0},\ldots,\bx_{m} \) that appear in \( F \).
\end{defn}

  Note that every monomial in \( F_{m} \) has weight \(m \) for $f\in k[[\bx_0]]$.

 The next lemma follows from  the definition of \( F_{m} \):
  
\begin{lem}
\label{appear}
\label{lem} Let \( f \) be a non-zero power series in \( k[[\bx_{0}]] \) of order \( 
\geq 1 \).
\begin{enumerate}
\item[(i)]
When char \( k \)= 0, 
a monomial \( \prod_{\ell=1}^r x_{0,j_{\ell}} \) appears  in \( f \) if and only if 
for every \( i_{\ell}\geq 0 \), 
the monomial \[ \prod_{\ell=1}^r x_{i_{\ell},j_{\ell}} \] appears in \( 
F_{m} \), where \( \sum_{\ell}i_{\ell}=m \).

Hence,  \( \ord F_{m}=\ord f \), and in particular \( F_{m}\neq 0 \) 
for every \( m \). 

\item[(ii)] 
For any characteristic,  a monomial $\displaystyle\prod_{j=1}^N x_{0,j}^{e_j}$ appears in $f$ if and only if 
for every \( i_{\ell}\geq 0 \), 
the monomial 
$$\prod_{j=1}^N x_{i_j,j}^{e_j}$$
appears  in $F_m$, where $m=\sum_j e_j i_j$.

\end{enumerate}

\end{lem}  

\begin{pf}
  The statement of ``if'' part follows immediately from the 
  definition of  \( F_{m} \) for both (i) and (ii). 
  Now assume that 
  \( g= \prod_{\ell=1}^r x_{0,j_{\ell}} \) is a monomial in \( f \).
  By substituting \( \sum_{i\geq 0}x_{i, j}t^i \) 
  into $x_{0,j}$ in  this monomial, we obtain 
  \[ g(\sum_{i\geq 0}\bx_{i}t^i)=G_{0}+G_{1}t+G_{2}t^2+\cdots. \]
  Therefore, \( G_{m} \) is the sum of the monomials of the form \( \prod_{\ell=1}^r x_{i_{\ell},j_{\ell}} \)
  with \( i_{\ell}\geq 0 \) and \( \sum_{\ell}i_{\ell}=m \).
  If the characteristic of \( k \) is zero, the coefficients of each such
  monomial is nonzero.
  And each monomial \( \prod_{\ell=1}^r x_{i_{\ell},j_{\ell}} \) in 
  \(G_{m}\) is not canceled by the contribution from the other monomials of \( f 
  \),
  because the collection \( (j_{1},..,j_{\ell},..,j_{r}) \)   
  assigns the source monomial \( \prod_{\ell=1}^r x_{0,j_{\ell}} \).
  This shows the statement of ``only if'' part of (i).
 For the proof of only if part of (ii), let $g=\prod_j x_{0,j}^{e_j}$ and define $G_i$  in the same way as in the previous discussion.
 Then, the monomial $\prod_j x_{i_j,j}^{e_j}$ appears with coefficient 1 in $G_m$ for $m=\sum_j e_j i_j$.
 Therefore, the coefficient of $\prod_j x_{i_j,j}^{e_j}$ in $F_m$ is the same as the coefficient of 
 $\prod_j x_{0,j}^{e_j}$ in $f$.
\end{pf}

\begin{rem} The statement (i) of Lemma \ref{appear} does not hold for positive characteristic case. For example, let $p>0$ be the characteristic of the base field $k$ and 
 \( f=x_{0,1}^p \in k[[x_{0,1}]]\).
 Then \( F_{m}=x_{i,1}^p \) for \( 
m=pi \) and \( F_{m}=0 \) for \( m\not\equiv 0 \) (mod \( p \)).
\end{rem}

As we saw in the previous section, 
Corollary \ref{smooth} follows immediately from Proposition \ref{sm}.
But here we give another proof of Corollary  \ref{smooth} for an algebraically closed base field, since we think  that it   gives some useful insight into jet schemes.
\vskip.5truecm
\noindent
[{\it Proof of Corollary \ref{smooth}}]
We may assume that  \( (X,0)\subset (\AAA, 0) \) is a closed 
subscheme with a singularity at $0$, where \( N \) is the embedding dimension of \( (X,0) \).
Then every element \( f  \in I_{X}  \) has order greater than 1. 
By this, every element  \( F_{i}\) of the defining ideal \( I_{X_{m}} \) 
of  \( X_{m} \) in \( (\AAA)_{m} \) 
 has order greater than 1. 
Here, note that \( I_{X_{m}}\neq 0 \), since \( I_{X}\neq 0 \) and \( 
F_{0}=f \) for \( f\in I_{X} \).
Therefore the Jacobian matrix of \( I_{X_{m}} \) is 
 the zero matrix at the trivial \( m 
\)-jet \( 0_{m} \in X_{m} \) at \( 0 \), which shows that \( 0_{m} \) is a singular 
point in \( X_{m} \) for every \( m \). \( \Box \)

\section{Proofs of theorems \ref{flat}, \ref{positive}}

\begin{say}
\label{note}
For the proof of the  theorems, 
we fix the notation as follows:
Let \( (X,0)\subset (\AAA, 0) \) be a singularity of 
embedding dimension \( N \).
Let \( 0\leq m<m' \), 
 \( R_{m}=k[[\bx_{0}]][\bx_{1},\ldots,\bx_{m}] \), \( I\subset R_{m} \) 
the defining ideal of \( X_{m} \) in \( (\AAA)_{m} \), \( 
R_{m'}=k[[\bx_{0}]][\bx_{1},..,\bx_{m},..\bx_{m'}] \) and \( I'\subset R_{m'} \) 
the defining ideal of \( X_{m'} \) in \( (\AAA)_{m'} \).
Let \( M \) be the maximal ideal of \( R_{m}\) generated by \( 
\bx_{0},\ldots,\bx_{m} \).
\end{say}

\begin{lem}
\label{notation}
Under the notation as in \ref{note}, if there is an element \( F\in I' 
\cap MR_{m'} \) such that  
\( F\not\in MI'+IR_{m'} \),
then the truncation morphism \(\psi_{m',m}: X_{m'}\to X_{m} \) 
 is not flat.
\end{lem}

\begin{pf} 
The truncation morphism \(\psi_{m',m}: X_{m'}\to X_{m} \) corresponds to 
the canonical  
ring homomorphism \( R_{m}/I\to R_{m'}/I' \).
The non-flatness follows from the non-injectivity of the 
canonical homomorphism:
\[ M/I\otimes_{R_{m}/I}R_{m'}/I'\to R_{m'}/I'. \]
Since we have an isomorphism of the first module
\[  
M/I\otimes_{R_{m}/I}R_{m'}/I'\simeq MR_{m'} /(MI'+IR_{m'}),\]
the existence of an element 
\( F\in I'\cap MR_{m'} \) such that \( F\not\in MI'+IR_{m'} \) gives 
the non-injectivity.
\end{pf}

[{\it Proof of Theorem \ref{flat}}]
Assume that the base field  \( k \) is algebraically closed and of characteristic zero and \( (X, 0) \) is 
a singular point of a scheme $X$ of finite type over $k$. 
Then we will deduce that every truncation morphism \(\psi_{m',m}: X_{m'}\to X_{m} 
\) 
\( (m'>m\geq 0) \) is not flat.
For this, it is sufficient to prove that  \(\psi_{m',m}: \widehat{X_{m'}}\to \widehat{X_{m}} 
\) 
\( (m'>m\geq 0) \) is not flat by Corollary \ref{reduction}.
So we may assume that $X$ is a closed subscheme of $\AAA$ with  the embedding dimension $N$.
Let $I_X$ be the defining ideal of $X$ in $\AAA$.
We use the notation of \ref{note}.
Let \( f  \) be an element in \( I_{X} \) with the minimal order \( d 
\).
Note that \( d\geq 2 \), as \( N \) is the embedding dimension.
Then, by Lemma \ref{lem}, (i), \( F_{m+1} \) is not zero and 
presented as 
\[ F_{m+1}=g_{1}(\bx_{0})x_{m+1, 1}+\cdots+g_{N}(\bx_{0})x_{m+1,N}+
g'(\bx_{0},\ldots,\bx_{m}), \]
where \( \ord 
F_{m+1}=d  \) and 
some of \( g_{i} \)'s are not zero.
We should note that \( \ord g_{i}=d-1 \) for all non-zero \( g_{i} \)'s.
As \( \ord g_{i}\geq 1\), for every \( i \) and \( \ord g'\geq 1 \), the element \( F_{m+1} \) is in \( MR_{m'} \).
It is  clear that \( F_{m+1}\in I' \).
On the other hand, 
as \( \ord I=\ord I'=d \), it follows that $\ord MI'\geq d+1$ and the initial term of an element 
$IR_{m'}$ of order $d$ is the initial term of an element of $I$.
Hence, 
the initial term of  an element in 
 \(  MI'+IR_{m'} \)  of order \( d \) should be the initial term of an element of $I$,
 therefore it should be a  polynomial in \( \bx_{0},\ldots, \bx_{m} \).
 However, the initial term of \( F_{m+1} \) is not of this form, 
 which     implies 
  \( F_{m+1}\not\in MI'+IR_{m'}\).
 By Lemma \ref{notation},  the non-flatness of \( 
  \psi_{m',m}:X_{m'}\to X_{m} \) follows for every pair \( (m, m') \) 
  with \( 0\leq m<m' \).
\( \Box \)

\begin{exmp}
\label{ex}
  The condition char\( k \)=0 is necessary for Theorem \ref{flat}. 
  Indeed, 
  there are counter examples for Theorem \ref{flat} in case of 
  positive characteristic.
  For example,
  let \( X \) be a scheme defined by \( x_{0,1}^p \) in \( 
  {\Bbb A}_{k}^1=\spec k[x_{0,1}] \) over a field \( k \) of characteristic p.
   Let $r$ be an integer with $0<r < p$
  Then, for any positive integer $q$, we have 
   \[ X_{pq+r}=\spec k[x_{0,1}, x_{1,1},..., x_{pq+r, 1}]/(x_{0,1}^p,..., x_{q,1}^p) \] and
  \[ X_{pq}=\spec  k[x_{0,1},   x_{1,1},.., x_{pq, 1}] / (x_{0,1}^p,..., x_{q,1}^p).  \]
  It is clear that \( X_{pq+r} \) is flat over \( X_{pq} \), while \( X \) is 
  singular.
\end{exmp}

\vskip.5truecm
\noindent
[{\it Proof of Theorem \ref{positive}}] 
As in the proof of the previous theorem, we will show the non-flatness of the truncation morphisms,
if $(X,0) $ is singular.
As $X$ is reduced, some fiber of the truncation morphism $\psi_{m',m}:X_{m'}\to X_m$ has 
dimension $\leq (m'-m)\dim (X,0)$ for a small affine neighborhood $X$ of $0$, if $\psi_{m',m}$ is flat.
(If $X$ is of equi-dimensional, then the fiber has dimension $\dim (X,0)$.)
Hence, if   $\psi_{m',m}$ is flat, by Corollaries \ref{equal}, \ref{reduction}, the dimension of the fiber over 
a closed point in $(\pi_{m}^{\widehat X})^{-1}(0)$ by the morphism
 $\widehat{\psi_{m',m}}:\widehat{X_{m'}}\to \widehat{X_m}$ is $\leq (m'-m)\dim (X,0)$.
With remarking this fact and Corollary \ref{reduction}, 
 we may assume that $X$ is a singular closed subscheme of 
$\AAA$ for the embedding dimension $N$ of $(X,0)$.

First  assume \( m'<d(m+1) \). 
  Note that for every \( g\in I_{X} \),  
 \[  \overline{G}_{i}=G_{i}({\bf 0},\ldots,{\bf 0},\bx_{m+1},\ldots,\bx_{i})=0 \] 
  for \( i<d(m+1) \).
  This is because every monomial in $G_i$ has a factor $x_{\ell, j}$ with $\ell\leq m$,
   since the weight of  \( G_{i} \) is \( i\) \( (<d(m+1))\) and 
   \(\ord G_{i}\geq d \).
  Let \( 0_{m} \) be the trivial \( m \)-jet at \( 0 \).
  As \( \psi_{m', m}^{-1}(0_{m}) \) is defined in \( \bA^{(m'-m)N} \) 
  by the ideal generated by \( \overline{G}_{i} \)'s with \( i\leq m' \) for $g\in I_X$,
  it follows that \[ \psi_{m', m}^{-1}(0_{m})\simeq \bA^{N(m'-m)}, \] 
  which is a fiber of dimension \( N(m'-m)> (m'-m)\dim (X,0) \).
 Therefore, $\psi_{m',m}$ is not flat,
 because otherwise the fiber dimension would be $(m'-m)\dim(X,0)$ as we saw before.

  Therefore, we may assume that \( m'\geq d(m+1) \), where $d=\ord I_X$.
  Let $f\in I_X$ have the order $ d$.
  Let $\prod_j x_{0,j}^{e_j}$ be a monomial with the minimal degree in $f$.
  Then, $\sum_j e_j=d$ and therefore $e_j\leq d$ for every $j$.
  Let $e$ be one of  non-zero $e_j$'s.
  By the assumption $m'\geq d(m+1)$, there is a positive integer $i$ such that $m\leq ie <m'$.
  Let $s$ be  minimal among such $i$'s.
   Then \( F_{se}\in I' \) is clear and also we have \( F_{se}\in MR_{m'} \) under the  notation of  \ref{note}.
  Indeed, if  a monomial \( \prod_{\ell=1}^u x_{i_{\ell}, j_{\ell}} 
  \) of \( F_{se} \) has a factor \( x_{i_{\ell}, j_{\ell}} \) with \( i_{\ell}\geq m+1 
  \),  let this \( i_{\ell} \) be \( i_{1} \).
  Then \( i_{1}\geq m+1 > (s-1)e  \).
  By this,
  \[ \sum_{\ell\neq 1}i_{\ell}< se-(s-1)e=e\leq d\leq u. \]
  Therefore, there is at least one \( \ell  \) such that \( i_{\ell}\leq 1\leq m \).
  Hence every monomial of \( F_{se} \) is contained in \( MR_{m'} \).
  Now let $e=e_1$.
  As 
   \[ \prod_j x_{0,j}^{e_j} \] 
  is a monomial  of  \( f \) of the minimal order \( d \), by Lemma \ref{appear},  
  \[ x_{1,s}^e\prod_{j\neq 1}x_{0, j}^{e_j} \]
  is a monomial of \( F_{se} \).
  Therefore, \( \ord F_{se}=d \).
  This monomial does not appear in any element of \( MI'+IR_{m'} \).
  Indeed,  \( \ord MI'\geq d+1 \) and the initial term of 
  an element of \( IR_{m'} \) of order \( d \)
  must be the initial term of an element of \( I \), because of \( 
  \ord I = d \). 
  Therefore, every initial monomial of  an element of \( IR_{m'} \) 
  of order \( d \) is of the 
  form
  \[ \prod_{\ell}x_{i_{\ell},j_{\ell}}, \ \ \  (\sum_{\ell}i_{\ell}\leq m),\]
  since \( I \) is generated by \( F_{i} \)'s with \( i\leq m \) for \( 
  f\in I_{X} \).
  As \(  x_{1,s}^e\prod_{j\neq 1}x_{0, j}^{e_j} \) is not of 
  this form, we obtain \( F_{se}\not\in IR_{m'}+MI' \).
 By this and Lemma \ref{notation}, it follows that
  \( X_{m'}\to X_{m}  \) is not flat for \( m'>m> 0 \).

 \begin{rem}
 In the proof of Theorem \ref{positive}, we used the condition $m\geq 1$.
 It is not clear if the same statement as in the Theorem \ref{positive}  follows for $m=0$
 in positive characteristic case, 
 i.e., If the base field is of positive characteristic, $X$ is reduced and 
 $\pi_{m'}=\psi_{m',0}:X_{m'}\to X$ is flat for some $m'>0$, then is $X$  non-singular?
 
  But in particular, if $m'=1$, it holds true.  
  This is seen as follows:
  For an affine scheme $X$ of finite type over $k$, the fiber of a point $x\in X$ by the projection $\pi_1:X_1\to X$ is the Zariski tangent space 
 of the point. 
 Therefore $\dim \pi_1^{-1}(x)=\ed (X,x)$.
 If $(X, 0)$ is singular and reduced, $\dim \pi_1^{-1}(x)> \dim(X,0)$, while there are points in a small neighborhood of $0$  such that the fiber dimension is $\dim (X,0)$. 
 Hence, $\pi_1$ is not flat.

 \end{rem}

\makeatletter \renewcommand{\@biblabel}[1]{\hfill#1.}\makeatother

\end{document}